\documentclass[12pt]{paper}

\usepackage{amssymb,amsfonts,amsthm,amsmath}

\usepackage{latexsym,multicol,graphicx}

\usepackage{graphicx, verbatim}
\usepackage{xcolor}
\usepackage{verbatim}

\usepackage[a4paper, hmargin=3cm, vmargin={3.9cm, 3.9cm}]{geometry}

\newtheorem{theorem}{Theorem}
\newtheorem*{theorem*}{Theorem}
\newtheorem{lemma}[theorem]{Lemma}

\renewcommand{\phi}{\varphi}

\renewcommand{\emptyset}{\varnothing}

\newcommand{\la}{\lambda}

\renewcommand{\le}{\leqslant}
\renewcommand{\ge}{\geqslant}

\newcommand{\bR}{\mathbb R}
\newcommand{\bC}{\mathbb C}

\newcommand{\bD}{\mathbb D}
\newcommand{\bT}{\mathbb T}

\numberwithin{equation}{section}
\numberwithin{theorem}{section}

\begin{document}

\title{A piece of Victor Katsnelson's mathematical biography}

\author{Mikhail Sodin}

\maketitle

\begin{abstract}
We give an overview of several works of Victor Katsnelson published in 1965--1970, and
pertaining to the complex and harmonic analysis and the spectral theory.
\end{abstract}

\section*{A preamble}
As a mathematician, Victor Katsnelson was raised within a fine school
of function theory and functional analysis, which was blossoming in
Kharkov starting the second half of 1930s. He studied in the
Kharkov State University in 1960-1965. Among his teachers were Naum Akhiezer,
Boris Levin, Vladimir Marchenko. That time he became acquainted with
Vladimir Matsaev whom Victor often mentions as one of his teachers.
In 1965 Katsnelson graduated with the master degree, Boris Levin supervised his
master thesis. Since then and till 1990,
he teaches at the Department of Mathematics and Mechanics of the Kharkov
State University. In 1967 he defends the PhD Thesis ``Convergence and Summability
of Series in Root Vectors of Some Classes of Non-Selfadjoint Operators'' also written
under Boris Levin guidance. Until he left Kharkov in the early 1990s, Katsnelson
remained an active participant of the Kharkov function theory seminar run on Thursdays
by Boris Levin and Iossif Ostrovskii. His talks, remarks and questions were always
interesting and witty.

Already in the 1960s Victor established himself among the colleagues as one of the finest
Kharkov mathematicians of his generation, if not the finest one. Nevertheless, he was not
appointed as a professor and was never allowed to travel abroad.

Most of Katsnelson's work pertain to the spectral theory of functions and operators.
I will touch only a handful of his results, mostly published in 1965--1970,
that is, at the very beginning of his mathematical career. A big portion of
his works written in Kharkov appeared in the local journal
``Function Theory, Functional Analysis and Their Applications'' and were never
translated in English. Today, the volumes of this journal are available at~{\tt http://dspace.univer.kharkov.ua/handle/123456789/43}.

In this occasion, let me mention two wonderful books carefully written
by Katsnelson~\cite{Kats-book1, Kats-book2}. They exist only as manuscripts, and curiously, both have ``Part~I'' in their titles, though, as far as I know, no continuations appeared.
In both books mathematics interlaces with interesting historical comments.
Last but not least, let me also mention an extensive survey of Issai Schur's works in
analysis written jointly by Dym and Katsnelson~\cite{Dym-Kats}.

\subsubsection*{Acknowledgements}
Alex Eremenko, Gennadiy Feldman,
Alexander Kheifets, Yuri Lyubarskii, Sasha Sodin, and Peter Yuditsii read a preliminary version
of this text and made several comments which I took into account. I thank all of them.

\section{A Paley-Wiener-type theorem}

The paper~\cite{Kats-Paley-Wiener} was, probably, the first published work of Katsnelson. Therein, he studied the following question raised by
Boris Levin. Given a convex compact set $K\subset\bC$ with the boundary $\Gamma=\partial K$, let  $\mathsf{L}^2(\Gamma)$ be the $L^2$-space of function
on $\Gamma$ with respect to the Lebesgue length measure. {\em How to characterize entire functions $F$ represented by the Laplace integral
\begin{equation}\label{eq:Laplace}
F(z) = \frac1{2\pi{\rm i}}\, \int_\Gamma  f(w) e^{w z}\, {\rm d} w,
\end{equation}
with $f\in \mathsf{L}^2(\Gamma)$}?

In the case when $K=\Gamma$ is an interval, the answer is provided by the classical Paley-Wiener theorem. In this case, it is convenient to assume
that $\Gamma\subset {\rm i}\bR$. Then we can rewrite~\eqref{eq:Laplace} as follows
\[
F(z) = \frac1{2\pi{\rm i}}\, \int_{{\rm i}a}^{{\rm i}b} f(w) e^{w z}\, {\rm d}w = \frac1{2\pi}\, \int_a^b \phi(t) e^{{\rm i}tz}\, {\rm d}t\,,
\qquad \phi\in\mathsf{L^2}(a, b)\,,
\]
and, by the Paley-Wiener theorem, a necessary and sufficient condition for this representation
with some $a<b$ is that $F$ is an entire function of exponential type (EFET, for short) and
$F\in\mathsf{L^2(\bR)}$.

Now, assume that the convex compact $K$ is not an interval, that is, is a closure of its interiour, and put $\Omega_K = \overline{\bC}\setminus K$. Note that the Laplace transform
of $F$  coincides with the Cauchy integral of $f$:
\[
\int_0^\infty F(z) e^{-\la z}\, {\rm d}z = \frac1{2\pi{\rm i}}\, \int_\Gamma \frac{f(w)}{\la-w}\, {\rm d} w\,.
\]
The RHS is analytic in $\Omega_K$, vanishes at infinity, and belongs to the Smirnov space $\mathsf{E}^2(\Omega_K)$,
which can be defined, for instance, as the closure in $\mathsf{L}^2(\Gamma)$ of
analytic functions in $\Omega_K$, continuous up to the boundary, and vanishing at infinity. Thus, Levin's question can be reformulated as follows: {\em Given a convex compact set $K$ with non-empty interiour, find a complete normed space $\mathsf{B}_K$
of EFET such that the Laplace integral $\mathcal L$ defined in~\eqref{eq:Laplace}
gives a bounded bijection $\mathsf{E}^2(\Omega_K) \stackrel{\mathcal L}\to \mathsf{B}_K$}.
Note that the representation~\eqref{eq:Laplace}
yields that all functions $F\in \mathsf{B}_K$ have the growth bound
\[
|F(re^{{\rm i}\theta})| \le C(\Gamma)\, \| f \|_{\mathsf{L}^2(\Gamma)}
\, \exp\bigl[ \max_{w\in K } \operatorname{Re}(we^{{\rm i}\theta}) r \bigr]
=  C(\Gamma)\, \| f \|_{\mathsf{L}^2(\Gamma)}\, e^{h_K(-\theta) r},
\]
where $h_K(\theta)$ is the supporting function of $K$.

The first result in that direction is due to Levin himself who considered in~\cite[Appendix~I, Section~3]{Levin-Distribution} the case when $K$ is a convex polygon and noticed that
in this case the answer is a straightforward consequence of the classical version of the Paley-Wiener theorem.  Then, M.~K.~Liht~\cite{Liht} considered the case when $K$ is a disk
centered at the origin  and of radius $h$. He showed that in this case one can take
$\mathsf{B}_K$ being a Bargmann-Fock-type space, which consists of entire functions $F$ satisfying
\[
\int_0^\infty\, \int_{-\pi}^\pi |F(re^{{\rm i}\theta})|^2\, e^{-2hr}\, \sqrt{r}\, {\rm d}r\, {\rm d}\theta <\infty\,.
\]

The starting point of Katsnelson's work~\cite{Kats-Paley-Wiener} was a remark that a more accurate version of the Liht argument yields an isometry
\[
\int_\Gamma |f|^2\, |{\rm d}w| = \int_0^\infty\, \int_{-\pi}^\pi |F(re^{{\rm i}\theta})|^2\, e^{-2hr}\, \rho(hr)\, {\rm d}r\, {\rm d}\theta\,,
\]
where
\[
\rho(r) = 2r \int_0^\infty \frac{e^{-2tr}}{\sqrt{2t+t^2}}\, {\rm d}t\,.
\]
Then, Katsnelson proves that representation~\eqref{eq:Laplace} yields a uniform bound
\[
\sup_{|\theta|\le\pi}\, \int_0^\infty |F(re^{{\rm i}\theta})|^2 e^{-2h_K(-\theta)r}\, {\rm d}r  \le C(\Gamma)\, \| f \|_{\mathsf{L}^2(\Gamma)}^2\,.
\]
The proof is based on the following lemma close in the spirit to known estimates due to
Gabriel and Carlson.
\begin{lemma}\label{lemma_Gabriel-Carlson}
Suppose that $K$ is a convex compact set, $\Gamma=\partial K$,
$\Omega_K=\overline\bC\setminus K$, and $f\in\mathsf{E}^2(\Omega_K)$.
Then, for any supporting line $\ell$ to $\Gamma$,
\[
\int_\ell |f|^2\, |{\rm d} w| \le C(\Gamma)\, \int_\Gamma |f|^2\, |{\rm d} w|\,.
\]
The constant on the RHS does not depend on $\ell$ and $f$.
\end{lemma}

One can modify Levin's question replacing the space $\mathsf E^2(\Omega_K)$ by
another space of functions analytic in $\Omega_K$. If functions in that space do
not have boundary values on $\Gamma$, then one needs to replace the integral over $\Gamma$
on the RHS of~\eqref{eq:Laplace} by the contour integral
\[
\frac1{2\pi{\rm i}}\, \int_\gamma f(w) e^{wz}\, {\rm d} w,
\]
where $\gamma$ is a simple closed contour in $\Omega_K$, which contains $K$ in its
interiour. This integral is called the Borel transform of $f$. It acts on the Taylor coefficients as follows:
\[
f(w) = \sum_{n\ge 0} \frac{a_n}{w^{n+1}} \mapsto F(z) = \sum_{n\ge 0} \frac{a_n}{n!}\, z^n.
\]

One of possible modifications of Levin's question
was considered in~\cite{Kats-Paley-Wiener}. Katsnelson introduces the weight
\[
\rho_\Gamma (w) = \frac1{|w-a_1(w)| + |w-a_2(w)|}\,,
\]
where $a_j(w)$, $j=1, 2$, are supporting points  for the line supporting to $\Gamma$ that passes through $w\in\Omega_K$ (the weight $\rho_\Gamma (w)$ is not defined when $w$ belongs to the supporting line to $\Gamma$ that has a common segment with $\Gamma$). The last result proven in~\cite{Kats-Paley-Wiener} is a curious isometry
\[
\int_0^\infty\, \int_{-\pi}^\pi |F(re^{{\rm i}\theta})|^2\, e^{-2h_K(-\theta)r}\,
{\rm d}r\, {\rm d}\theta = \frac1{2\pi} \iint_{\Omega_K} |f(w)|^2 \rho_\Gamma (w)\,
{\rm d}\sigma(w)\,,
\]
where $\sigma$ is the Lebesgue area measure.

\medskip Works of Liht and Katsnelson had follow-ups. In~\cite{Lyubarskii},
Lyubarskii extended Liht's theorem to convex compact sets $K$ with smooth boundary.
The decisive word was said by Lutsenko and Yulmukhametov. In~\cite{Luts-Yulm2}
they proved that the Laplace integral $\mathcal L$ defines an isomorphism\footnote{Understood as an isomorphism
between Banach spaces} between
$\mathsf{E}^2(\Omega_K)$ and a space of EFET such that
\[
\int_0^\infty\, \int_{-\pi}^\pi |F(re^{{\rm i}\theta})|^2\,
\frac{{\rm d}r\, {\rm d}\Delta(\theta)}{K(re^{{\rm i}\theta})}\,,
\]
where
\[
K(z) = \| e^{w z} \|^2_{\mathsf{E}^2(\Omega_K)}
= \int_\Gamma e^{2\operatorname{Re}(w z)}\, |{\rm d}w|\,,
\]
and ${\rm d}\Delta (\theta) = (h''(-\theta) + h(-\theta)){\rm d}\theta$ (understood as
a distribution). One of the novelties in their work is the fact that the identity map
provides an isomorphism between the Smirnov space
$\mathsf{E}^2(\Omega_K)$ and the space of analytic
functions in $\Omega_K$ vanishing at infinity, with finite Dirichlet-type
integral
\[
\iint_{\Omega_K} |f'(w)|^2 \operatorname{dist}(w, \Gamma)\, {\rm d}\sigma (w) < \infty.
\]
The proof of that fact relies on Lemma~\ref{lemma_Gabriel-Carlson}.
We mention that Yulmukhametov together with his pupils and collaborators
proved several other non-trivial results related to Levin's question
(see, for instance,~\cite{Luts-Yulm1, Isaev-Yulm}) and that Lindholm~\cite{Lindholm}
extended the Lutsenko-Yulmukhametov theorem to analytic functions of several complex variables.

\section{Riesz bases of eigenvectors of non-selfadjoint operators}

One of the central questions in the spectral theory is the expansion in eigenfunctions
(more generally, of root vectors) of non-selfadjoint operators. It
originates in the theories of ordinary and partial differential equations and of integral equations. In the middle of the 1960s the corresponding completeness problem was already
understood relatively well, first of all, due to the pioneering works by Keldysh and
Matsaev. A portion of their works can be found
in the classical Gohberg-Krein book~\cite{Go-Kr}, another portion became available
later in~\cite{Keldysh} and in~\cite{Matsaev-Mogulskii1, Matsaev-Mogulskii2}.
The situation with convergence of the series of eigenfunctions was understood much less clear. Though a few results, due to Glazman, Mukminov,
and Markus, were known (all of them were summarized in~\cite[Chapter~VI]{Go-Kr}), no general
methods existed until in~\cite{Kats-root_vectors} Katsnelson discovered
a novel approach to the Riesz basis property of eigenfunctions of arbitrary contractions
and dissipative operators. His approach is based on a deep result of Carleson pertaining to the interpolation by bounded analytic functions in the unit disk.

We start with some definitions. First, we remind the notion of Riesz basis of a system
of subspaces $(\mathsf X_k)$ of a Hilbert space $\mathsf H$. The details can be found in~\cite[Chapter~VI]{Go-Kr}. In the case when all subspaces $(\mathsf X_k)$ are one-dimensional,
this notion reduces to the usual notion of the Riesz basis of vectors in $\mathsf H$.

Let $(\mathsf X_k)$ be a collection of linear subspaces of $\mathsf H$, and $\mathsf X$ be the closure of
their linear span. {\em The subspaces $(\mathsf X_k)$ form a basis in $\mathsf X$} if any vector
$x\in\mathsf X$ has a unique decomposition into a convergent series
$ x = \sum_k x_k, \quad x_k\in \mathsf X_k $.
To simplify notation, we assume that the linear span of the subspaces $(\mathsf X_k)$ is dense in $\mathsf H$, i.e., that $\mathsf X=\mathsf H$.
Let $\mathsf P_k$ be projectors on $\mathsf X_k$. Then the system $(\mathsf X_k)$
forms a basis if and only if $\mathsf P_k \mathsf P_j = \delta_{kj}\mathsf P_k$,
and $ \sup_n \bigl\| \sum_{k=1}^n \mathsf P_k \bigr\| < \infty $.
The subspaces $(\mathsf X_k)$ form {\em an orthogonal basis} if all $\mathsf P_k$s are orthogonal projectors, that is, for any $x$,
$\| x\|^2 = \| \mathsf P_k x \|^2 + \| (\mathsf I - \mathsf P_k)x\|^2$.
The subspaces $(\mathsf X_k)$ form
{\em a Riesz basis} if there exists an invertible operator $\mathsf A$ from $\mathsf H$ onto
$\mathsf H$ such that
subspaces $(\mathsf A \mathsf X_k)$ form an orthogonal basis.
Gelfand's theorem~\cite[Chapter~VI, \$~5]{Go-Kr} says that {\em a basis of subspaces
$(\mathsf X_k)$ is a Riesz basis if and only if it remains a basis after any permutation of
its elements}.

In~\cite{Kats-root_vectors} Katsnelson studies the question when the collection of root subspaces of a non-selfadjoint operator is a Riesz basis in the closure of its linear span.
He considers two general classes of non-selfadjoint operators, contractions and dissipative operators.
A linear operator $\mathsf T$ on a Hilbert space $\mathsf H$
is called {\em a contraction} if $\| \mathsf T \|\le 1$. Given an eigenvalue $\la\in\overline\bD$,
the linear space
\[
\mathsf X(\la) = \bigcup_{n\ge 1} \operatorname{Ker}\bigl[ (\mathsf T-\la \mathsf I)^n \bigr]
\]
is called {\em the root subspace} corresponding to $\la$.
The eigenvalue $\la$ {\em has finite order} if there exists a positive integer
$m$ such that
\[
\mathsf X(\la) =
\bigcup_{1\le n \le m} \operatorname{Ker}\bigl[ (\mathsf T-\la \mathsf I)^n \bigr]
= \operatorname{Ker}\bigl[ (\mathsf T-\la \mathsf I)^m \bigr].
\]
The least value $m$ is called {\em the order $m(\la)$ of the eigenvalue $\la$}.
The following theorem is the main result of~\cite{Kats-root_vectors}.

\begin{theorem}\label{thm-root_vectors}
Let $(\la_k)$ be some eigenvalues of a contraction $\mathsf T$, let $(\mathsf X_k)$
be the corresponding root subspaces, and let $(m_k)$ be the orders of $(\la_k)$.
Suppose that
\begin{equation}\label{eq_separation}
\inf_j\, \prod_{k\ne j} \Bigl| \frac{\la_k-\la_j}{1-\la_j \bar\la_k} \Bigr|^{m_j m_k}
\ge \delta >0\,,
\end{equation}
Then the system of root subspaces
$(\mathsf X_k)$ forms a Riesz basis in the closure of its linear span.
\end{theorem}
In~\cite{Kats-root_vectors} Katsnelson only sketches the proof of this result,
some details can be found in Nikolskii's survey paper~\cite[\S~3]{Nik}.
Here are the main steps of the proof.

First, Katsnelson observes that when in the assumptions of Theorem~\ref{thm-root_vectors}
only those $\la_k$ that lie in the open unit disk matter, while the unitary part of
the operator $\mathsf T$ can be discarded. He also assumes that
the linear span of the root subspaces $\mathsf X(\la_k)$ is dense in
$\mathsf H$ (otherwise, he considers
the restriction of $\mathsf T$ on the closure of this linear span). Keeping in mind
Gelfand's theorem, it suffices to find projectors
$\mathsf P_j\colon \mathsf H\to \mathsf X(\la_j)$ such that
$\mathsf P_j \mathsf X(\la_k)=\{0\}$ for $j\ne k$, and
\[
\sup_J \bigl\|\, \sum_{j\in J} \mathsf P_j\, \bigr\| < \infty\,,
\]
where the supremum is taken over all finite subsets $J$ of the set of all indices $j$.

Fix a finite set $J$. Suppose that we succeeded to find an analytic in the unit disk
function $f_J$ such that
$ f_J(\la_j)=1 $ for $j\in J$,
$ f_J^{(\nu)}(\la_j)=0$ for $j\notin J$ and $0\le \nu \le m(\la_j)-1$,
and $ \sup_\bD |f_J| \le M $,
with a constant $M$ independent of $J$.
Suppose momentarily that the function $f_J$ is analytic on a neighbourhood of the closed unit disk (that is, that the set $(\la_j)$ is finite). Then, by the F.~Riesz operator calculus~\cite[Chapter~IX]{RS},
$f_J(\mathsf T)$ is well-defined and, by von
Neumann's theorem~\cite[Section~153, Theorem~A]{RS},  $\| f_J(\mathsf T) \| \le M$.

At the next step, Katsnelson again uses a piece of the F.~Riesz operator calculus. The projectors $\mathsf P_j$ can be defined by the contour integrals
\[
\mathsf P_j
= - \frac1{2\pi {\rm i}}\, \int_{C_j} (\mathsf T-\zeta \mathsf I)^{-1}\, {\rm d}\zeta\,,
\]
where $C_j$ is a circumference of a small radius which separates the point $\la_j$
from the rest of the spectrum and traversed counterclockwise. Whence,
\[
f_J(\mathsf T)
= - \frac1{2\pi {\rm i}}\, \int_{\bT} f_J(\zeta)(\mathsf T-\zeta \mathsf I)^{-1}\,
{\rm d}\zeta = \sum_{j\in J} \mathsf P_j\,,
\]
and therefore,
\[
\bigl\|\, \sum_{j\in J} \mathsf P_j\, \bigr\| = \| f_J(\mathsf T) \|
= \sup_{\bD} |f_J| \le M\,.
\]

To get rid of the assumption that the function $f_J$ is analytic on the neighbourhood of
the closed unit disk, Katsnelson applies a classical result due to Pick and Schur,
which says that
{\em given a bounded analytic function $f$ in the unit disk and given a finite set of points
$\Lambda\subset\bD$, there exists a rational function $R$ which interpolates $f$ at $\Lambda$,
that is, $R(\la)=f(\la)$, $\la\in\Lambda$, and $\max_{\bD} |R| = \sup_\bD |f|$}, see, for instance,~\cite[Corollary~IV.1.8]{Garnett}.

 At the final step, Katsnelson deduces the existence of the analytic function $f_J$ with the properties as above from Carleson's ``$0-1$-interpolation theorem'', which, in turn, was the main step in his solution to the corona problem~\cite[Theorem~2]{Carleson}.
\hfill $\Box$

\medskip
This chain of arguments discovered in~\cite{Kats-root_vectors}
had a significant impact on works of many mathematicians, notably
from the Saint Petersburg school, cf.
Nikolskii-Pavlov~\cite{Nik-Pavlov1, Nik-Pavlov2} (apparently, Nikolskii and Pavlov
rediscovered some of Katsnelson's results), Treil~\cite{Treil1, Treil2}, Vasyunin~\cite{Vasyunin}, see also~\cite[Lectures~IX and X]{Nik-book}.

\medskip Katsnelson also notes that
condition~\eqref{eq_separation} in Theorem~\ref{thm-root_vectors} cannot be weakened.
Given a sequence $(\la_k)\subset\bD$ satisfying the Blaschke condition $\sum_k (1-|\la_k|)<\infty$ and such that
\[
\inf_j\, \prod_{k\ne j} \Bigl| \frac{\la_k-\la_j}{1-\la_j \bar\la_k} \Bigr|^{m_j m_k} = 0,
\]
he brings a simple construction (the idea of which, according to~\cite{Kats-root_vectors},
is due to Matsaev) of a contraction $\mathsf T$ such that

\smallskip\noindent (i) $(\la_k)$ are simple eigenvalues of $\mathsf T$ and
the whole spectrum of $\mathsf T$ coincides with $(\la_k)$, and

\smallskip\noindent (ii) the eigenvalues of $\mathsf T$ are complete in $\mathsf H$
but are not uniformly minimal\footnote{
A system of vectors $\{x_n\}$ in the Hilbert space $\mathsf H$ is called {\em uniformly minimal}
if there exists $\delta>0$ such that for all $n$ the distance between $x_n$ and the linear span of $\{x_k\colon k\ne n\}$ is at least $\delta$.}
.

\smallskip\noindent Furthermore, the operator $\mathsf I-\mathsf T^* \mathsf T$ is one-dimensional.

\medskip Among other results brought in~\cite{Kats-root_vectors}, there is a version
of Theorem~\ref{thm-root_vectors} for dissipative operators, i.e., the operators
$\mathsf A$ such
that $\operatorname{Im} \langle \mathsf Ax, x \rangle \ge 0 $, for any
$x$ in the domain of $\mathsf A$.
This version is reduced to Theorem~\ref{thm-root_vectors} by an application of the Caley transform $\mathsf A\mapsto (\mathsf A-{\rm i}\mathsf I)(\mathsf A+{\rm i}\mathsf I)^{-1}$.

\section{Series of simple fractions}

Let $\mathsf{C_0(\bR)}$ be the Banach space of complex-valued continuous functions on $\bR$,
tending to zero at infinity, equipped with the uniform norm $\| f \| = \sup_\bR |f|$.
Fix finite subsets in the upper and lower half-planes
$\{z_k\}_{1\le k \le n}\subset\bC_+$ and
$\{w_k\}_{1\le k \le m}\subset\bC_-$ and denote by
$\mathsf{E_+} = \mathsf{E_+(w_1, \ldots , w_m)}$,
$\mathsf{E_-} = \mathsf{E_-(z_1, \ldots , z_n)}$ the subspaces in
$\mathsf{C_0(\bR)}$
generated by the simple fractions $\{ 1/(t-w_k)\colon 1\le k \le m \}$ and
$\{ 1/(t-z_k)\colon 1\le k \le n \}$. The functions in
$\mathsf{E_+}$ are analytic on $\bC_+$,
the functions in $\mathsf{E_-}$ are analytic on $\bC_-$.
Furthermore, $\mathsf{E_+} \cap \mathsf{E_-} = \emptyset$ and the sum
$\mathsf E = \mathsf{E_+} + \mathsf{E_-}$ is a direct one,
i.e., for any function $f\in\mathsf E$, there exists a unique decomposition $f=f_+ + f_-$ with $f_\pm \in \mathsf{E_\pm}$.
In~\cite{Kats-simple_fractions} Katsneslon estimates the norms of the projectors
$\mathsf{P_\pm} = \mathsf{P_\pm (z_1, \ldots , z_n; w_1, \ldots , w_m)}$ from $\mathsf E$
onto the corresponding subspace $\mathsf{E_\pm}$. The main result of that work
is the following theorem:

\begin{theorem}\label{thm-simple_fraction}
There exists a positive numerical constant $C$ such that
\[
\| \mathsf{P_\pm} \| \le C\min (m, n)\, (m+n)\,.
\]
\end{theorem}
The main point in this theorem is that the upper bound it gives does not depend on the
positions of $z_k$s and $w_k$s.

Note that in the space $\mathsf{L^2(\bR)}$, by one of the versions of the Paley-Wiener theorem, the functions analytic in the upper and lower half-planes are orthogonal to each other, which makes the corresponding projectors orthogonal. In view of this remark, it is quite natural that the proof of Theorem~\cite{Kats-simple_fractions} uses the Fourier transform. The proof
is nice and not too long and the reader can find its details in~\cite{Kats-simple_fractions}.
Bochtejn and Katsnelson bring in~\cite{Boch-Kats} a counterpart of Theorem~\ref{thm-simple_fraction} for the unit circle $\bT$ instead of the real line $\bR$.

\medskip Likely, the upper bound given in Theorem~\ref{thm-simple_fraction} is not sharp. In~\cite{Kats-simple_fractions}, Katsnelson conjectures that, for $m=n$, the sharp upper bound should be $C\, \log (n+1)$, and, as far as we know, this conjecture remains open till today.
As a supporting evidence towards this conjecture, he brings the following result.

\begin{theorem}\label{thm:simplest}
Given two finite sets of points $\{z_k\}_{1\le k\le n}\subset\bC_+$ and
$\{w_k\}_{1\le k\le n}\subset\bC_-$, consider the functions
\[
g_+(t) = \sum_{k=1}^n \frac1{t-w_k}\,,
\quad g_-(t) = \sum_{k=1}^n \frac1{t-z_k}\,,
\]
and let $ g(z) = g_+(z) + g_-(z) $.
Then
\[
\max_\bR |g_\pm| \le C\log (n+1)  \cdot \max_\bR |g|\,.
\]
with a positive numerical constant $C$.
\end{theorem}

A simple example shows that the order of growth of the RHS cannot be improved. Put
\[
g_+(t) = \sum_{k=1}^n \frac1{t+{\rm i} k}\,,
\quad g_-(t) = \sum_{k=1}^n \frac1{t-{\rm i} k}\,.
\]
Then
\[
g(t) = g_+(t) + g_-(t) = \sum_{k=1}^n \frac{2t}{t^2+k^2}\,,
\]
and
\[
\max_\bR |g| \le \int_0^\infty \frac{2t}{t^2+x^2}\, {\rm d}x = \pi,
\]
while
\[
\max_\bR |g_-| = |g_-(0)| =\sum_{k=1}^n \frac1k = \log n + O(1)\,.
\]

The proof of Theorem~\ref{thm:simplest} is short and elegant (and accessible to undergraduate
students). As a byproduct of that proof he obtains
\begin{theorem}\label{thm-log_derivative}
Let $P$ be a polynomial of degree $n\ge 2$ such that
\[
\sup_\bR \Bigl| \frac{P'}{P} \Bigr| \le M\,.
\]
Then $P$ has no zeroes in the strip
\[
|\operatorname{Im} z| \le \frac{c}{M \log n}\,,
\]
where $c$ is a positive numerical constant.
\end{theorem}
Slightly earlier a similar estimate was obtained by Gelfond~\cite{Gelfond}. Gelfond's proof was rather different (and more involved).
The question about the size of the strip around the real axis free of
zeroes of $P$ was raised by Gorin in~\cite{Gorin62}, first results in that
direction were obtained by him and then by Nikolaev in~\cite{Nikolaev}.
The final word in this question was said by Danchenko, who proved
in~\cite{Danchenko} that under assumption of Theorem~\ref{thm-log_derivative},
the polynomial $P$ has no zeroes in the strip
\[
|\operatorname{Im} z| \le \frac{c}{M}\cdot \frac{\log\log n}{\log n}\,,
\]
and that the order of decay of the RHS cannot be improved.

\medskip Twenty five years later, Katsnelson returned in~\cite{Kats-Tumarkin1, Kats-Tumarkin2} to the linear spans of simple fraction but from a different point of view.
That time his work was motivated by Potapov's results on factorization of $J$-contractive
matrix functions.

Let $m$ be the Lebesgue measure on the unit circle $\bT$, and
$w\colon \bT\to [0, \infty]$ be an $m$-integrable weight, satisfying the Szeg\H{o}
condition
\begin{equation}\label{eq:Szego}
\int_\bT \log w\, {\rm d}m > -\infty.
\end{equation}
By $\mathsf{PCH^2(w)}$ Katsnelson denotes the Hilbert space of functions
$f$ analytic on $\bC\setminus\bT$ and satisfying the following conditions

\smallskip\noindent (a) The restriction of $f$ onto $\bD_+$ and $\bD_-$ belongs to
the Smirnov class, i.e., $\log_+|f|$ has positive harmonic majorants both in
$\bD_+$ and $\bD_-$.

\smallskip\noindent (b) The boundary values of $f\big|_{\bD_+}$ and $f\big|_{\bD_-}$
coincide $m$-a.e. on $\bT$, that is,
\[
\lim_{r\uparrow 1} f(rt) = \lim_{r\downarrow 1} f(rt)\ (\, =: f(t)\, )
\quad m-{\rm a.e.\ on\ } \bT.
\]
(Conditions (a) and (b) together provide the so called pseudocontinuation property
of the function $f$.)

\smallskip\noindent (c)
\[
\| f \|_w^2 = \int_\bT |f|^2 w\, {\rm d} m < \infty\,.
\]

\medskip Note that whenever $w^{-1}\in\mathsf{L}^1(m)$ the space $\mathsf{PCH^2(w)}$
is trivial, i.e., contains only constant functions $f$. Indeed, convergence of the integrals
\[
\int_{\bT} |f|^2 w\, {\rm d}m < \infty, \quad
\int_{\bT} \frac{{\rm d}m}{w} < \infty
\]
yields that $f\in\mathsf{L}^1(m)$, and then, by a version of the removable singularity theorem that goes back to Carleman, the function $f$ is entire, and since it is bounded, by Lioville's theorem it is a constant function.

\medskip Given a set of points $S\subset \bD_+\cup\bD_-$ satisfying the Blaschke condition
\begin{equation}\label{eq:Blaschke}
\sum_{\la\in S\cap\bD_+} (1-|\la|) < \infty\,,
\quad \sum_{\la\in S\cap\bD_-} (1-|\la|^{-1}) < \infty\,,
\end{equation}
denote by $\mathsf{R(S; w)}$ the closure of the linear span of the simple fractions
$\{(t-\la)^{-1}\}_{\la\in S}$ together with the constant functions in the space
$\mathsf{L^2(w)}$. Let $S=S_1 \supset S_2 \supset \, \ldots \, $ be a chain of sets
such that $\bigcap_n S_n = \{\emptyset\}$.

The starting point of Katsnelson's work~\cite{Kats-Tumarkin1} is the inclusion
$ \bigcap_n \mathsf{R(S_n; w)} \subset \mathsf{PCH^2(w)} $, which follows
from classical results of Tumarkin~\cite{Tumarkin1},
see also~\cite{Tumarkin2, Tumarkin3}. Katsnelson observes that this inclusion might be a strict one, that is, generally speaking, not every function in
the space $\mathsf{PCH^2(w)}$ can be approximated in $\mathsf{L^2(w)}$ by a sequence of
functions $r_n\in\mathsf{R(S_n; w)}$.
The main result of~\cite{Kats-Tumarkin1} is the following approximation theorem.

\begin{theorem}\label{thm:tumarkin}
For any non-negative $m$-integrable function $w$ on $\bT$ satisfying the Szeg\H{o}
condition~\eqref{eq:Szego}, there exists a set $S\subset\bD_+\cup\bD_-$ satisfying the Blaschke condition~\eqref{eq:Blaschke} such that
$ \bigcap_n \mathsf{R(S_n; w)} = \mathsf{PCH^2(w)} $.
\end{theorem}

In~\cite{Kats-Tumarkin2} Katsnelson extends this result to a more general approximation scheme by simple fractions with poles at a given table of points in $\mathbb C\setminus\bT$.
In~\cite{Khefa} Kheifets used Katsnelson's construction to answer a question raised by Sarason.

\section{Spectral radius of hermitian elements in Banach algebras and
the Bernstein inequality}

One of the most important properties of EFET is the classical Bernstein inequality, which
states that {\em if $F$ is an entire function of exponential type $\sigma$, then
\[
\sup_{\bR} |F'| \le \sigma\, \sup_{\bR} |F|\,,
\]
and the equality sign attains if and only if $F(z) = c_1 \cos\sigma z + c_2\sin\sigma z$.}
Different proofs, deep extensions, and various applications of the Bernstein inequality
can be found in the books~\cite{Akhiezer, Levin-Distribution, Levin-Lectures} and in the
survey paper~\cite{Gorin}. Interestingly, Bernstein's inequality is also closely related to the theory of Banach algebras.

An element $a$ of a Banach algebra $\mathsf A$ is called {\em hermitian} if
$\| e^{{\rm i}at}\| =1$ for every $t\in\bR$. For instance, hermitian elements of the algebra of all bounded operators in a Hilbert space are self-adjoint operators. Another, more special, example
is the differentiation operator $  \mathsf D = \tfrac1{\rm i}\, \tfrac{\rm d}{{\rm d}x} $ considered in various Banach spaces of EFET equipped with some translation-invariant norm,
in which case the exponent $ e^{{\rm i}\mathsf D t} $ is realized by the translation.
It is well-known that the operator norm of a self-adjoint operator in a Hilbert space coincides with its spectral radius. Making use of the Bernstein inequality, Katsnelson proved in~\cite{Kats-Bernstein} the following result, which, independently (and more or less simultaneously), was also found by~Browder~\cite{Browder} and Sinclair~\cite{Sinclair}.

\begin{theorem}\label{thm-spectral_radius}
For every hermitian element in a Banach algebra, the norm coincides with the spectral radius. \end{theorem}

Moreover, as both Katsnelson and Browder observed, this result is {\em equivalent} to the Bernstein inequality, that is, {\em the latter follows from the former, applied to the differentiation operator $\mathsf D$ in the Bernstein space $\mathsf{B_\sigma}$ of EFET at most $\sigma$ bounded on the real axis and equipped with the uniform norm}.

The proof of Theorem~\ref{thm-spectral_radius} is short and elegant:
Let $a$ be a hermitian element in a Banach algebra $\mathsf A$. Take an arbitrary linear functional $\phi\in\mathsf A^*$ with the unit norm, and consider the EFET
\[
F(z) \stackrel{\rm def}= \phi\bigl( e^{{\rm i}a z} \bigr) = \sum_{n\ge 0} \phi(a^n)\, \frac{z^n}{n!}\,.
\]
Applying, first, the formula, which expresses the exponential type of an entire function
via its Taylor coefficients, then a crude estimate of $n!$, and then Gelfand's formula
for the spectral radius, we estimate the exponential type of $F$:
\begin{multline*}
\sigma_F = \frac1{e}\, \limsup_{n\to\infty} n\Bigl( \frac{|\phi(a^n)|}{n!} \Bigr)^{1/n} \\
\le \frac1{e}\, \limsup_{n\to\infty} n\Bigl( \frac{\| a^n \|}{n!} \Bigr)^{1/n}
= \limsup_{n\to\infty} \| a^n \|^{1/n} = \rho (a)\,,
\end{multline*}
where $\rho(a)$ denotes the spectral radius of $a$. Since the element $a$ is hermitian,
we have $|F(x)| = |\phi\bigl( e^{{\rm i}ax} \bigr)| \le \| e^{{\rm i}ax} \| = 1$, whence,
by the Bernstein inequality,
\[
|\phi(a)| = |F'(0)| \le  \sigma_F \sup_\bR |F| \le \rho (a)\,,
\]
and then, by the Hahn-Banach theorem, $\| a \| \le \rho(a)$. This completes the proof of Theorem~\ref{thm-spectral_radius} since the converse inequality
$\rho (a) \le \| a \|$ is obvious. \hfill $\Box$

\medskip The proof of the result which goes in the opposite direction is also quite simple.
Let $\mathsf D$ be the differentiation operator in the Bernstein space
$\mathsf{B_\sigma}$ of EFET at most $\sigma$ bounded on $\bR$. As we have already mentioned, the exponential function $e^{{\rm i}\mathsf Dt}$ acts on $\mathsf{B_\sigma}$
as the translation by $t$, so $\mathsf D$ is a hermitian operator in $\mathsf{B_\sigma}$. To evaluate the spectral radius of $\mathsf D$,
we need to estimate from above the norms $\| \mathsf D^n \|$, that is,
$\sup_\bR |F^{(n)}|$, $F\in\mathsf{B_\sigma}$. By Cauchy's estimate for the derivatives of analytic functions, combined with the bound $|F(x+w)|\le e^{\sigma |w|}\, \sup_\bR |F|$
valid for any $F\in\mathsf{B_\sigma}$, we obtain
$ |F^{(n)}(x)| \le n!\, r^{-n} e^{\sigma r}\| F \| $ for any $r>0$ and any $x\in\bR$. Optimising the RHS, we get $ |F^{(n)}(x)| \le n! \exp[n-n\log n + n\log\sigma]\, \| F \|$,
that is, $ \| \mathsf D^n \| \le n! \exp[ n - n\log n + n\log\sigma ] $, and finally,
$\rho (\mathsf D) = \lim \| \mathsf D^n \|^{1/n} = \sigma$.
Thus, for any function $F\in\mathsf{B_\sigma}$
and any $x\in\bR$, we have
\[
|F'(x)| = | (\mathsf DF)(x) | \le \| \mathsf D \|\cdot \|F\|
= \rho(\mathsf D)\| F \| \le \sigma\| F \|,
\]
proving the Bernstein inequality. \hfill $\Box$

\medskip In this context, it is also worth mentioning that a bit later Bonsall and Crabb~\cite{Bonsall-Crabb} found a simple direct proof of Theorem~\ref{thm-spectral_radius}, which yields another proof of the Bernstein inequality. Their proof is based on the
following lemma, which is a simple exercise on the functional calculus in Banach algebras:
\begin{lemma}\label{lemma:arcsin}
Let $a$ be a hermitian element in a Banach algebra with $\rho (a) < \pi/2$. Then,
$a=\arcsin(\sin a)$.
\end{lemma}
Now, Theorem~\ref{thm-spectral_radius} follows almost immediately. Proving Theorem~\ref{thm-spectral_radius}, it suffices, assuming that $a$ is an arbitrary
hermitian element with $\rho(a)<\pi/2$, to show that $\| a \| \le \pi/2 $.
Let $c_n$ be the Taylor coefficient of the function $z\mapsto \arcsin z$, $|z| \le 1$.
The values $c_n$ are positive and their sum equals $\arcsin (1) = \pi/2$.
By Lemma~\ref{lemma:arcsin}, $ \| a \| \le \sum_{n\ge 1} c_n \| \sin a \|^n$. Since the element $a$ is hermitian, $\| \sin a\|\le 1$, and therefore,
$\| a \| \le \sum_{n\ge 1} c_n = \pi/2$. \hfill $\Box$

\centerline{* \quad * \quad *}

In the reference list, referring to the papers in Russian published in journals translated from cover to cover, we mention only the translations. Today, the original Russian versions of these papers can be found at the Math-Net.Ru site ({\tt http://www.mathnet.ru}).

\bigskip
\medskip

\noindent
School of Mathematics, Tel Aviv University, Tel Aviv 69978, Israel
\newline email:{\tt\ sodin@tauex.tau.ac.il}

\end{document}